\documentclass[12pt]{article}
\usepackage{amssymb}

\newenvironment{rem}{\bigskip \noindent \bf Remark \rm}{\bigskip}
\newenvironment{example}{\bigskip \noindent \bf Example \rm}{\bigskip}
\newenvironment{proof}{\bigskip \noindent \bf Proof: \rm}{\bigskip}

\newenvironment{enumthm}{
\begin{enumerate}}{\end{enumerate}}
\newtheorem{thm}{Theorem}
\newtheorem{cor}[thm]{Corollary}
\newcommand{\qed}{\begin{flushright} \vspace{-1pc} $\square$
\end{flushright}}
\newfont{\knot}{virtual scaled 1000}
\begin{document}
\begin{center}
{\Large\bf On Alexander-Conway Polynomials \\[1ex]
for Virtual Knots and Links} \\[3ex]
\large J\"org Sawollek\footnote{Fachbereich Mathematik, Universit\"at
Dortmund, 44221 Dortmund, Germany \\ {\em E-mail:\/}
sawollek@math.uni-dortmund.de \\ {\em WWW:\/}
http://www.mathematik.uni-dortmund.de/lsv/sawollek} \\[1ex]
December 21, 1999 (revised: January 3, 2001)
\end{center}
\vspace{4ex}

\begin{abstract}
A polynomial invariant of virtual links, arising from an invariant of
links in thickened surfaces introduced by Jaeger, Kauffman, and Saleur,
is defined and its properties are investigated. Examples are given that
the invariant can detect chirality and even non-invertibility of virtual
knots and links. Furthermore, it is shown that the polynomial satisfies
a Conway-type skein relation -- in contrast to the Alexander polynomial
derived from the virtual link group.  \\[1ex]
\end{abstract}

\section*{Introduction}
In \cite{kauf99} Kauffman defines an extension of classical knot
diagrams to virtual knot diagrams, motivated by Gauss codes on the one
hand and knots in thickened surfaces on the other hand. Several
classical knot invariants can be generalized to the virtual theory
without much effort, e.g., the knot group and derived invariants such as
the Alexander polynomial, the bracket and Jones polynomials, and
Vassiliev invariants (which can be introduced in different ways, see
\cite{kauf99} and \cite{gpv}).

The present paper deals with a polynomial invariant that is derived from
an invariant of links in thickened surfaces introduced by Jaeger,
Kauffman, and Saleur in \cite{jks}. The determinant formulation of the
polynomial immediately generalizes to virtual link diagrams. It is a
Laurent polynomial in two variables with integral coefficients that
vanishes on the class of classical link diagrams but gives non-trivial
information for diagrams that represent non-classical virtual links.
Especially, examples can be given that the invariant is sensitive with
respect to changes of orientation of a virtual knot. Furthermore, the
polynomial fulfills a Conway-type skein relation in one variable and
thus it is denoted by the term Conway polynomial.

In the same way as in the classical case, the one-variable Alexander
polynomial of a virtual link can be derived from the virtual link group,
but the skein-relation for (a normalized version of) the classical
Alexander polynomial cannot be extended to the class of virtual links.
Therefore, this Alexander polynomial is different from the Conway
polynomial mentioned above in a non-trivial way -- in contrast to the
classical case (and certain generalizations to links 3-manifolds, see
for example \cite{kai}, Theorem 5.2.11).

This paper is organized as follows. After, in Section \ref{virtknots}, a
short introduction into the field has been given, the determinant
formulation of the Conway polynomial for virtual links is described in
Section \ref{invariant}. Some properties of the polynomial are deduced,
especially, that the invariant fulfills a Conway-type skein relation,
and several example calculations are given. Then, in Section \ref{alex},
the Alexander invariants derived from the link group, namely, Alexander
matrix, Alexander ideal, and Alexander polynomial, are defined and it is
shown that the Alexander polynomial does not fulfill any linear skein
relation. Finally, in Section \ref{conrem}, general problems in
extending certain invariants of classical links to the virtual category
are described and the direction of further investigations is indicated.

\section{Virtual Knots and Links}
\label{virtknots}

\begin{figure}[htb]
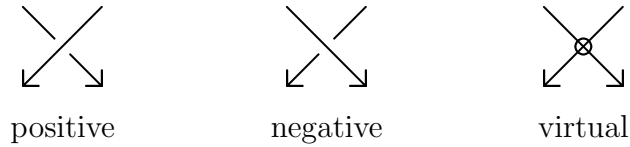

\vspace{-2ex}

\begin{center}
\begin{tabular}[t]{c@{\hspace{5em}}c@{\hspace{5em}}c}
\\[1ex]
{\knot p} & {\knot q} & {\knot r} \\[1ex]
positive & negative & virtual
\end{tabular}
\end{center}
\vspace{-1ex}

\caption{Crossing types}
\label{crossings}
\end{figure}

In classical knot theory, knots and links in 3-dimensional space are
examined. As a main tool, projections of such links to an appropriate
plane are considered, namely, the so-called {\em link diagrams} (see,
for example, \cite{rol}, \cite{bur}, \cite{kaufbook}, \cite{kawa},
\cite{mura}, \cite{lick}). The idea of virtual knot theory is to
consider link diagrams where an additional crossing type is allowed.
Thus, for an oriented diagram, there are three types of crossings: the
classical positive or negative crossings and the virtual crossings (see
Fig. \ref{crossings}).
For a motivation to introduce virtual crossings resulting from examing
knots in thickened surfaces see \cite{kauf99}.

\begin{figure}[htb]
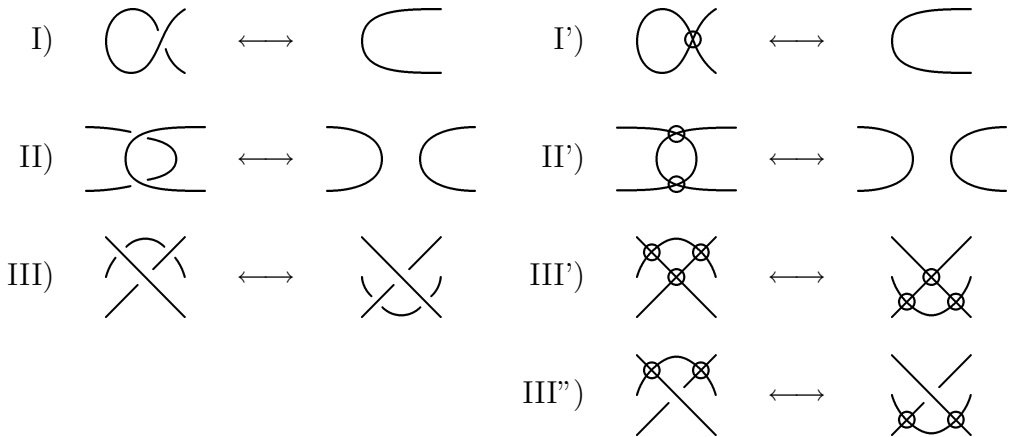

\begin{tabular}[t]{rccc}
I) & {\knot a} & $\longleftrightarrow$ & {\knot b} \\[2ex]
II) & {\knot d} & $\longleftrightarrow$ & {\knot e \hspace{0.5cm} f}
\\[2ex]
III) & {\knot g} & $\longleftrightarrow$ & {\knot h}
\end{tabular} \hspace{0em}
\begin{tabular}[t]{rccc}
I') & {\knot n} & $\longleftrightarrow$ & {\knot b} \\[2ex]
II') & {\knot m} & $\longleftrightarrow$ & {\knot e \hspace{0.5cm} f}
\\[2ex]
III') & {\knot k} & $\longleftrightarrow$ & {\knot l} \\[2ex]
III'') & {\knot i} & $\longleftrightarrow$ & {\knot j}
\end{tabular}
\vspace{2ex}

\caption{Extended Reidemeister moves}
\label{reidemeister}
\end{figure}

\begin{figure}[htb]
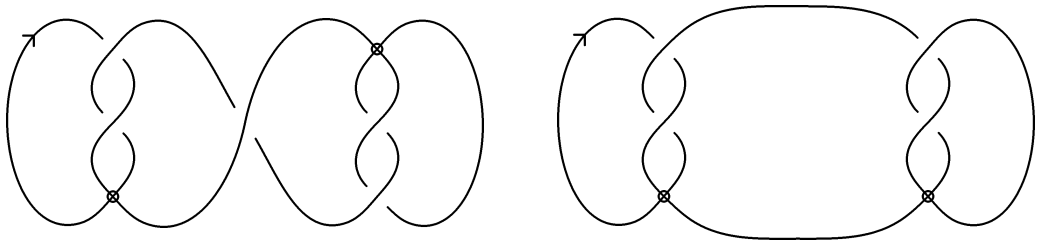

\begin{center}
\knot z \hspace{1cm} A
\end{center}
\caption{Two different virtual knots}
\label{flype}
\end{figure}

Technically, a {\em virtual link diagram\/} is an oriented 4-valent
planar graph embedded in the plane with appropriate orientations of
edges and additional crossing information at each vertex as depicted in
Fig. \ref{crossings}. Denote the set of virtual link diagrams by $\cal
VD$. Two diagrams $D, D' \in {\cal VD}$ are called {\em equivalent\/} if
one can be transformed into the other by a finite sequence of extended
Reidemeister moves (see Fig. \ref{reidemeister})
combined with orientation preserving homeomorphisms of the plane to
itself. A {\em virtual link\/} is an equivalence class of virtual link
diagrams. This gives a purely combinatorial definition of virtual links
and, indeed, it is not allowed to perform all the modifications of a
virtual link diagram that are inspired by moving objects in 3-space as
in classical knot theory. For example, the virtual knots depicted in
Fig. \ref{flype} 
are not equivalent, which can be shown by calculating their Conway
polynomials defined in Section \ref{invariant}, though they could
easily be seen to be equivalent by performing a "flype" if there were
arbitrary classical crossings instead of the two virtual ones.

In \cite{kauf99} the Jones polynomial is extended to virtual links by
making use of the connection with the bracket polynomial. The latter one
can be defined for virtual link diagrams in the same way as for
classical link diagrams using the bracket skein relation, i.e., every
classical crossing is cut open in the two possible ways which results in
diagrams having only virtual crossings and therefore being equivalent to
crossing-free diagrams.

\begin{figure}[htb]
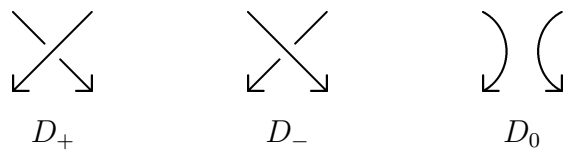

\begin{center}
\vspace{1ex}

\begin{tabular}[t]{c@{\hspace{5em}}c@{\hspace{5em}}c}
{\knot p} & {\knot q} & {\knot s} \\[1ex]
$D_+$ & $D_-$ & $D_0$
\end{tabular}
\end{center}
\vspace{-1ex}

\caption{Skein triple}
\label{skein}
\end{figure}

Extending other well-known link polynomials that can be defined via
skein relations, such as HOMFLY and Kauffman polynomials, is much more
difficult because it is in general not possible to get a diagram of a
trivial link from an arbitrary virtual link diagram by changing some of
the diagram's classical crossings. Therefore, it is a priori not clear
which basis can be chosen for the {\em skein module\/} (see \cite{hopr})
corresponding to a given skein relation. This problem also arises when
defining the Conway polynomial $\nabla_D(z) \in \mathbb Z[z]$ for a
diagram $D$ via the skein relation
\[ \nabla_{D_+}(z) - \nabla_{D_-}(z) = z \nabla_{D_0}(z) \]
where $(D_+, D_-, D_0)$ is a {\em skein triple\/}, i.e., the three
diagrams are identical except in a small disk where they differ as
depicted in Fig. \ref{skein}.

In the following two sections two different extensions of the classical
one-variable Alexander-Conway polynomial are examined. The first one
satisfies a Conway-type skein relation in one variable and thus is
called Conway polynomial, but it is introduced without using skein
theory. The second one is the Alexander polynomial derived from the
virtual link group.

\section{A Polynomial Invariant of Virtual Links}
\label{invariant}

A general method to define invariants of virtual links that very often
does work is to apply the definition of an invariant for classical link
diagrams to virtual link diagrams by ignoring the virtual crossings.
This method is used in what follows to define an invariant for virtual
links which is essentialy identical to an invariant for links in
thickened surfaces that has been introduced in \cite{jks}.

Let $D$ be a virtual link diagram with $n \geq 1$ classical crossings
$c_1, \ldots, c_n$. Define
\[ M_+ := \left( \begin{array}{cc}
1-x & -y \\[1ex]
-xy^{-1} & 0
\end{array} \right) \qquad \mbox{and} \qquad M_- := \left(
\begin{array}{cc}
0 & -x^{-1}y \\[1ex]
-y^{-1} & 1-x^{-1}
\end{array} \right) \, . \]
For $i = 1, \ldots, n$, let $M_i := M_+$ if $c_i$ is positive and let
$M_i := M_-$ otherwise. Define the $2n \times 2n$ matrix $M$ as a block
matrix by $M := diag(M_1, \ldots, M_n)$.

\begin{figure}[htb]
\begin{center}
\begin{picture}(1,1.5)
\put(0,0){\knot o}
\put(-0.5,1.2){$i_l^-$}
\put(-0.5,-0.2){$i_l^+$}
\put(1.2,1.2){$i_r^-$}
\put(1.2,-0.2){$i_r^+$}
\end{picture}
\end{center}
\caption{Half-edges belonging to the vertex $v_i$}
\label{vertex}
\end{figure}
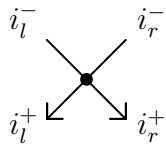

Furthermore, consider the graph belonging to the virtual link diagram
where the virtual crossings are ignored, i.e., the graph consists of $n$
vertices $v_1, \ldots, v_n$ corresponding to the classical crossings and
$2n$ edges corresponding to the arcs connecting two classical crossings
(the edges possibly intersect in virtual crossings). Subdivide each edge
into two half-edges and label the four half-edges belonging to the
vertex $v_i$ by $i_l^-$, $i_r^-$, $i_r^+$, $i_l^+$ as depicted in Fig.
\ref{vertex}.
A permutation of the set $\{ 1, \ldots, n \} \times \{ l, r \}$ is given
by the following assignment: $(i,a) \mapsto (j,b)$ if the half-edges
$i_a^+$ and $j_b^-$ belong to the same edge of the virtual diagram's
graph. Let $P$ denote the corresponding $2n \times 2n$ permutation
matrix where rows and columns are enumerated $(1,l)$, $(1,r)$, $(2,l)$,
$(2,r)$, \ldots, $(n,l)$, $(n,r)$.

Finally, define $Z_D(x,y) := (-1)^{w(D)} \det(M-P)$, where $w(D)$
denotes the {\em writhe\/} of $D$, i.e, the number of positive crossings
minus the number of negative crossings in $D$. (If $D$ has no classical
crossings then $Z_D(x,y)$ can be defined by $Z_D(x,y) := 0$, see Theorem
\ref{vanish}.)

\begin{thm}
$Z : {\cal VD} \to \mathbb Z[x^{\pm 1}, y^{\pm 1}]$ is an invariant of
virtual links up to multiplication by powers of $x^{\pm 1}$.
\end{thm}

\begin{proof}
The independence of the definition with respect to the ordering of the
classical crossings and the invariance of $\det(M-P)$ under Reidemeister
moves of type II and III follows exactly as in \cite{jks}, the only
difference being an exchange of the variable $-x$ for $x$. The behaviour
under Reidemeister moves of type I is depicted in Fig. \ref{reide1}.
\begin{figure}[htb]
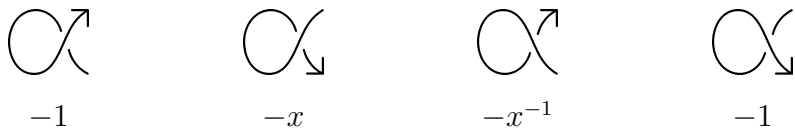

\begin{center}
\vspace{1ex}

\begin{tabular}[t]{c@{\hspace{5em}}c@{\hspace{5em}}c@{\hspace{5em}}c}
{\knot t} & {\knot u} & {\knot v} & {\knot w} \\[1ex]
$-1$ & $-x$ & $-x^{-1}$ & $-1$
\end{tabular}
\end{center}
\vspace{-2ex}

\caption{Behaviour of $Z$ under Reidemeister moves of type I}
\label{reide1}
\end{figure}
Since the changes of sign corresponding to Reidemeister moves of type I
are compensated by the factor $(-1)^{w(D)}$ in $Z_D(x,y)$ and since the
definition of $M$ and $P$ does not depend on the virtual crossings of
the diagram, the statement on $Z$ follows immediately.
\qed
\end{proof}

\begin{rem}
A normalization of the polynomial $Z_D(x,y)$ using rotation numbers, as
done in \cite{jks}, is not possible in the case of virtual link diagrams
because of Reidemeister move I'. But due to the change of the variable
"$x$" to "$-x$" in comparison with \cite{jks}, at least the sign of the
polynomial is determined.
\end{rem}

Define the {\em normalized polynomial\/} $\widetilde{Z}_D(x,y)$ as
follows. If $Z_D(x,y)$ is a non-vanishing polynomial and $N$ is the
lowest exponent in the variable $x$ then define $\widetilde{Z}_D(x,y) :=
x^{-N} Z_D(x,y)$. Otherwise let $\widetilde{Z}_D(x,y) := Z_D(x,y) = 0$.

\begin{cor}
$\widetilde{Z} : {\cal VD} \to \mathbb Z[x, y^{\pm 1}]$ is an invariant
of virtual links.
\end{cor}
\qed

\begin{example}
Let $D$ be the virtual knot diagram depicted in Fig. \ref{trefoil}.
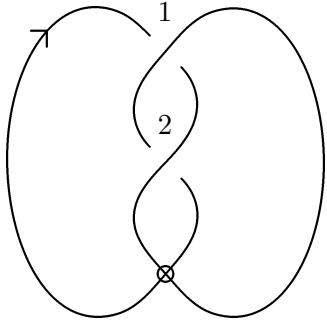
\begin{figure}[htb]
\begin{center}
\begin{picture}(4.5,4)
\put(0,-0.5){\knot y}
\put(2,3.7){\small 1}
\put(2,2.2){\small 2}
\end{picture}
\end{center}
\caption{A virtual trefoil}
\label{trefoil}
\end{figure}
The diagram arises from a diagram of the right-handed trefoil where a
positive crossing has been replaced by a virtual one. The corresponding
polynomial $Z_D(x,y)$ can be calculated by using the definition:
\[ (-1)^2 \det \! \left( \begin{array}{cccc}
1-x & -y & 0 & -1 \\[1ex]
-xy^{-1} & 0 & -1 & 0 \\[1ex]
-1 & 0 & 1-x & -y \\[1ex]
0 & -1 & -xy^{-1} & 0
\end{array} \right) = \, x^2+x^2y^{-1}+xy-xy^{-1}-y-1 \]
The normalized polynomial $\widetilde{Z}_D(x,y)$ is identical to the
polynomial $Z_D(x,y)$ in this case.
\end{example}

\begin{thm}
\label{vanish}
Let $D$, $D_1$, $D_2$ be virtual link diagrams and let $D_1 \sqcup D_2$
denote the disconnected sum of the diagrams $D_1$ and $D_2$. Then the
following hold.
\begin{enumthm}
\item $Z_D(x,y) = \widetilde{Z}_D(x,y) = 0$ if $D$ has no virtual
crossings
\item $Z_{D_1 \sqcup D_2}(x,y) = Z_{D_1}(x,y) Z_{D_2}(x,y)$,
$\widetilde{Z}_{D_1 \sqcup D_2}(x,y) = \widetilde{Z}_{D_1}(x,y)
\widetilde{Z}_{D_2}(x,y)$
\end{enumthm}
\end{thm}

\begin{proof}
Part a) follows as in \cite{jks}. Part b) is an immediate consequence of
the definition of the matrix $M-P$.
\qed
\end{proof}

\begin{rem}
For a connected sum $D_1 \# D_2$ of virtual link diagrams $D_1$ and
$D_2$, a formula of the form
\[ Z_{D_1 \# D_2}(x,y) = c Z_{D_1}(x,y) Z_{D_2}(x,y) \]
with a constant $c$ does not hold in general, in contrast to what could
be expected at first glance. For example, let $D_1$ be a diagram with
non-vanishing $Z$--polynomial and let $D_2$ be a diagram of the trivial
knot. Then the equation above obviously gives a contradiction.
\end{rem}

\begin{example}
The opposite of Theorem \ref{vanish} a) does not hold in general. A
counter-example is given in Fig. \ref{balance}.
\begin{figure}[htb]
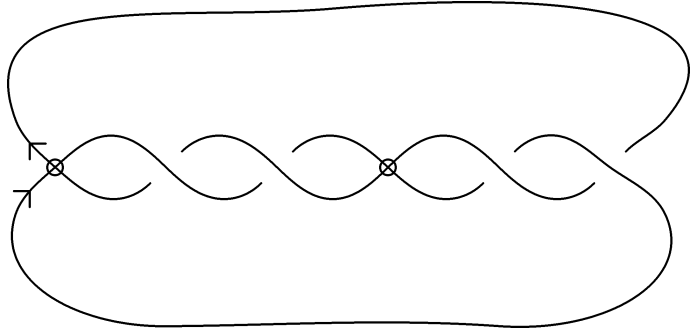

\begin{center}
\knot B
\end{center}
\caption{Non-classical virtual link with vanishing $Z$--polynomial}
\label{balance}
\end{figure}
The diagram represents a non-classical link since reversing the
orientation of one of the two components yields a virtual link diagram
with non-trivial $Z$--polynomial, whereas the $Z$--polynomial of the
diagram with the original choice of orientations vanishes.
\end{example}

\begin{thm}
\label{skeinrel}
Let ($D_+$, $D_-$, $D_0$) be a skein triple of virtual link diagrams.
Then the following skein relation holds:
\[ x^{-\frac{1}{2}} Z_{D_+}(x,y) - x^{\frac{1}{2}} Z_{D_-}(x,y) =
(x^{-\frac{1}{2}} - x^{\frac{1}{2}}) Z_{D_0}(x,y) \]
\end{thm}

\begin{proof}
The formula can easily be checked by verifying the corresponding
relation for every state $< \! v \, | \, f \! >$ in the state sum
model used in \cite{jks} to define the "partition function" $Z_D(x,y)$.
\qed
\end{proof}

Since a skein relation as in Theorem \ref{skeinrel} is fulfilled,
in the following $Z_D(x,y)$ and $\widetilde{Z}_D(x,y)$ will be called
{\em Conway polynomial\/} and {\em normalized Conway polynomial\/},
respectively. A more familiar version of the Conway skein relation can
be achieved by setting $x:=t^2$ and defining $Z_D'(t,y) := t^{-w(D)}
Z_D(t,y)$. Then, considering $w(D_0) = w(D_+)-1 = w(D_-)+1$, Theorem
\ref{skeinrel} immediately gives the following.

\begin{cor}
Let ($D_+$, $D_-$, $D_0$) be a skein triple of virtual link diagrams.
Then the following skein relation holds:
\[ Z'_{D_+}(t,y) - Z'_{D_-}(t,y) = (t^{-1} - t) Z'_{D_0}(t,y) \]
\end{cor}
\qed

\begin{rem}
Applying the skein relation of Theorem \ref{skeinrel} is, of course,
very helpful for calculating the Conway polynomial. For example, the
virtual link depicted in Fig. \ref{balance} can immediately be seen to
have vanishing Conway polynomial since changing an arbitrary crossing in
the diagram yields the diagram of a Hopf link and cutting open the same
crossing gives a diagram of the trivial knot. For the diagram that
arises from changing the orientation of one link component, the
corresponding skein tree has a branch that ends in a diagram with two
virtual and two classical crossings which cannot be simplified by any
changes of the classical crossings. Therefore, latter diagram has to be
calculated by using the definition of the Conway polynomial.
\end{rem}

\begin{thm}
\label{basis}
Let $D$ be a virtual link diagram. Then:
\begin{enumthm}
\item $Z_D(x,-x) = 0$
\item $Z_D(x,-1) = 0$
\item $Z_D(1,y)$ does not depend on the over-under information of the
diagram's classical crossings.
\end{enumthm}
\end{thm}

\begin{proof}
Part a) and b) follow from the fact that summing up the columns (rows)
of the determinant belonging to $Z_D(x,-x)$ ($Z_D(x,-1)$) gives the
trivial column (row) vector. Part c) is an immediate consequence of
setting $x=1$ in Theorem \ref{skeinrel} (or in the definition).
\qed
\end{proof}

\begin{rem}
By Theorem \ref{basis}, $Z_D(1,y) \in \mathbb Z[y^{\pm 1}]$ is an
invariant of virtual links that is invariant with respect to changes of
classical crossings, too. Thus $Z_D(1,y)$ can give information to
distinguish those "basic links" whose Conway polynomials must be
calculated by the determinant formulation instead of applying the skein
relation, i.e., generating elements of the virtual skein module related
to the Conway skein relation. For example, the value of $Z_D(1,y)$ for
the non-trivial diagram at the leaf of the skein tree that has been
mentioned in the previous remark is non-zero. Therefore, the
corresponding generating element is different from the trivial knot.
\end{rem}

For a virtual link diagram $D$, let $\overline{D}$ denote the {\em
mirror image\/} of $D$, i.e., the diagram that arises from $D$ by
changing the over-under information of every classical crossing. If $D$
and $\overline{D}$ are equivalent then $D$ is called {\em
amphicheiral\/} and otherwise {\em chiral}. Furthermore, let $D^*$ be
the {\em inverse\/} of $D$, i.e., the diagram with inverse orientation
for every component of $D$. If $D$ and $D^*$ are equivalent then $D$ is
called {\em invertible\/} and otherwise {\em non-invertible}.

In the case of classical links, it is known that Kauffman and HOMFLY
polynomials, and therefore the Conway polynomial too, are insensitive
with respect to invertibility of a knot or link (indeed, this is true
for any invariant derived from quantum groups, see \cite{kup},
\cite{lin}). Surprisingly, the Conway polynomial of a virtual knot can
be different from the polynomial of its inverse.

\begin{example}
The virtual knot with diagram $D$ that is depicted in Fig. \ref{notinv}
\begin{figure}[htb]
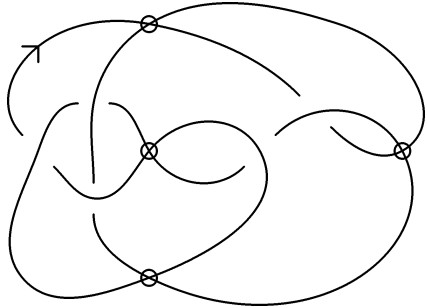

\begin{center}
\knot x
\end{center}
\caption{A chiral and non-invertible virtual knot}
\label{notinv}
\end{figure}
is chiral as well as non-invertible. Calculating normalized Conway
polynomials of $D$, $\overline{D}$, $D^*$, $\overline{D}^*$ gives the
following result:
\[ \widetilde{Z}_D(x,y) = (x-1)^2 (y+1) (x^2y^{-2}-1) \, , \;
\widetilde{Z}_{D^*}(x,y) = y \widetilde{Z}_D(x,y) \]
\[ \widetilde{Z}_{\overline{D}}(x,y) = - (x-1)^2 (y^2-1) (xy^{-1}+1) \,
, \; \widetilde{Z}_{\overline{D}^*}(x,y) = y^{-1} \widetilde{Z}_D(x,y)
\]
Thus the four diagrams under consideration represent pairwise different
virtual knots.
\end{example}

\section{The Alexander Polynomial Derived from the Link Group}
\label{alex}

In the same way as for classical links, the {\em (virtual) link group\/}
can be defined via the Wirtinger presentation of a virtual link diagram,
i.e., a group generator is assigned to each arc of the diagram and a
group relation can be read off at each classical crossing corresponding
to the rules shown in Fig. \ref{wirtinger}
\begin{figure}[htb]
\begin{center}
\begin{picture}(3,2.5)
\put(0,1){\vector(1,0){2}}
\put(1,0){\line(0,1){0.8}}
\put(1,1.2){\vector(0,1){0.8}}
\put(2.2,0.9){$b$}
\put(1.2,0.3){$a$}
\put(1.2,1.7){$c = bab^{-1}$}
\end{picture}
\hspace{8em}
\begin{picture}(3,2.5)
\put(2,1){\vector(-1,0){2}}
\put(1,0){\line(0,1){0.8}}
\put(1,1.2){\vector(0,1){0.8}}
\put(-0.4,0.9){$b$}
\put(1.2,0.3){$a$}
\put(1.2,1.7){$c = b^{-1}ab$}
\end{picture}
\end{center}
\caption{Wirtinger relations}
\label{wirtinger}
\end{figure}
(for details see \cite{kauf99}). The link group is an invariant of
virtual links.

The above definition of the link group is a purely combinatorial one and
it has, in general, nothing to do with the complement of a link in
3-space as in classical knot theory. Indeed, it can be shown that there
exists a virtual knot group which is not the fundamental group of any
3-manifold. Related work on virtual knot groups can be found in
\cite{kim00} and \cite{SiWi98}.

As explained in \cite{fox}, the {\em Alexander matrix\/} corresponding
to a group presentation can be calculated via a differential calculus.
For a presentation with $n$ generators and $m$ relations, it is a $m
\times n$ matrix with entries from the ring $\mathbb Z[t^{\pm 1}]$. The
{\em (first) Alexander ideal\/} ${\cal E}(D)$ of a diagram $D$ is
generated by the $(n-1) \times (n-1)$ minors of the Alexander matrix for
the corresponding link group and the {\em (first) Alexander
polynomial\/} $\Delta_D(t) \in \mathbb Z[t^{\pm 1}]$ is the greatest
common divisor of ${\cal E}(D)$. The Alexander ideal is an invariant of
virtual links and the Alexander polynomial is an invariant up to sign
and up to multiplication by powers of $t^{\pm 1}$.

\begin{rem}
In \cite{SiWi00} Alexander polynomials derived from an {\em extended
Alexander group\/} of a virtual link diagram are considered. Indeed, the
Conway polynomial $Z_D$ is related to the {\em 0th virtual Alexander
polynomial\/} of \cite{SiWi00}, see Remark 4.2
therein.
\end{rem}

In contrast to the classical Alexander polynomial, the Alexander
polynomial for virtual links does not fulfill any linear skein relation
as stated in the next theorem. Therefore it is crucially different from
the Conway polynomial defined in Section \ref{invariant}.

\begin{thm}
\label{noskein}
For any normalization $A_D(t)$ of the polynomial $\Delta_D(t)$, i.e.,
$A_D(t) = \varepsilon_D t^{n_D} \Delta_D(t)$ with some $\varepsilon_D
\in \{ -1, 1 \}$ and $n_D \in \mathbb Z$, the equation
\[ p_1(t) A_{D_+}(t) + p_2(t) A_{D_-}(t) + p_3(t) A_{D_0}(t) = 0 \;\;
\mbox{with } p_1(t), p_2(t), p_3(t) \in \mathbb Z[t^{\pm 1}] \]
has only the trivial solution $p_1(t) = p_2(t) = p_3(t) = 0$.
\end{thm}

\begin{proof} Assume the above skein relation has a non-trivial
solution. Consider the skein triples $(D_+, D_-, D_0)$ where $D$ is a
classical knot diagram with one crossing and $(D'_+, D'_-, D'_0)$ where
$D'$ is a standard diagram, with arbitrary orientations of the
components, of the Hopf link (the latter triple corresponding to an
arbitrary of the diagram's two crossings). Then the related Alexander
polynomials have values, up to normalization, as follows.
\[ \Delta_{D_+}(t) = \Delta_{D_-}(t) = \Delta_{D'_0}(t) = 1, \,
\Delta_{D_0}(t) = \Delta_{D'_-}(t) = 0, \, \Delta_{D'_+}(t) = t-1 \]
Inserting these values, each multiplied by a factor $\varepsilon t^n$,
into the skein relation immediately shows that, up to normalization,
$p_2(t) = p_1(t)$ and $p_3(t) = (t-1)p_1(t)$. Thus the skein relation is
equivalent to the classical Alexander skein relation:
\[ A_{D_+}(t) - A_{D_-}(t) = (t-1) A_{D_0}(t) \]
Finally, let $D''$ be the virtual link diagram arising from $D'$ by
changing an arbitrary classical crossing to a virtual crossing. Then
$\Delta_{D''_+}(t) = \Delta_{D''_-}(t) = t-1$ and $\Delta_{D''_0}(t) =
1$. Inserting these values, each multiplied by a factor $\varepsilon
t^n$, into the Alexander skein relation yields a contradiction.
\qed
\end{proof}

\begin{rem}
The proof of Theorem \ref{noskein} shows that the classical Alexander
skein relation is, up to a factor, the unique one that holds for the
Alexander polynomial of classical links, and this skein relation cannot
be extended to the class of virtual links.
\end{rem}

\begin{example}
The virtual knot with diagram $D$ which is depicted in Fig. \ref{notinv}
cannot be distinguished from its inverse by Alexander polynomials since
the Alexander ideals are identical:
\[ {\cal E}(D) = {\cal E}(D^*) = (2, t^2+t+1) \]
\end{example}

\begin{example}
\begin{figure}[htb]
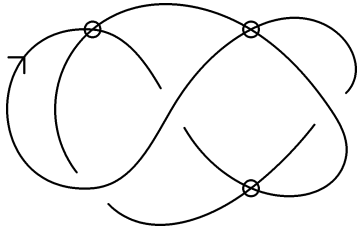

\begin{center}
\knot C
\end{center}
\caption{A virtual knot with trivial Jones polynomial}
\label{kaufex}
\end{figure}
The virtual knot with trivial knot group and trivial Jones polynomial
investigated in \cite{kauf99}, see Fig. \ref{kaufex}, has non-vanishing
normalized Conway polynomial
\[ (x-1) (x^2-y^2) (1+y^{-1}) \]
and therefore it does not represent a classical knot. This is a result
that could not be achieved with the means of \cite{kauf99}.
\end{example}

\section{Concluding Remarks}
\label{conrem}

As mentioned in \cite{jks}, the Alexander-Conway polynomial for links in
thickened surfaces can be defined in the multivariable case as well and
the construction analogously yields a multivariable Conway polynomial
for virtual links. Also the multivariable Alexander polynomial can be
derived from the virtual link group and it is an invariant that is
different from the Conway polynomial as has been seen above for the
one-variable case.

It is quite natural to consider generalizations of HOMFLY and Kauffman
polynomials to virtual links next. Besides the definition of these
polynomials via skein relations, there are several state models for
them, see \cite{CoRi}, \cite{jae89}, \cite{jae92}, \cite{jae93},
\cite{kauf90}, \cite{tur}, and also \cite{kaufbook} and citations
therein. But, when trying to generalize any of these approaches, one
meets with at least one of the following three obstacles. First, a
virtual link diagram in general cannot be unknotted by classical
crossing changes. Therefore it is necessary to find a basis for the
virtual skein module corresponding to the HOMFLY and Kauffman skein
relations, respectively, which is not yet known. Secondly, some models
make use of the signed graph corresponding to a black-and-white
colouring of the plane, but it is not clear how to handle virtual
crossings in a generalized model. And finally, state models mostly rely
on rotation numbers which are not invariant with respect to Reidemeister
moves of type I'.

It should be mentioned that missing invariance under Reidemeister moves
of type I and I' is a more serious problem when defining virtual link
invariants than missing invariance under Reidemeister moves of type I
when defining classical link invariants. For example, having in mind
Theorem \ref{noskein}, it is not possible to derive a linear skein
relation for the Alexander polynomial from the virtual link group. The
proof that is given by Hartley \cite{hart} for the classical case cannot
be extended to virtual link diagrams because it intrinsically uses a
normalization of the Alexander polynomial via rotation numbers during
the proof. Observe that in classical knot theory a result by B. Trace
\cite{tra} assures that a regular isotopy invariant can always be made
to an ambient isotopy invariant using writhe and rotation number of an
oriented link diagram.

A more general class of invariants that may be generalized to virtual
links are {\em quantum link invariants}. In \cite{kauf99} it is
described how to define such invariants via $R$-matrices and the
abstract tensor approach, see also \cite{kaufbook}. Again, the
invariance under Reidemeister moves of types I and I' causes
difficulties and thus the virtual quantum link invariants under
consideration are only invariant with respect to {\em virtual regular
isotopy\/}, i.e., the two troublemaking moves are avoided.

The Conway polynomial for virtual links that has been introduced in the
present paper also arises from a quantum link invariant (see
\cite{KaSa}) but it is not difficult to get control over its behaviour
with respect to Reidemeister moves of types I and I'. It is left to
further investigations which (quantum link) invariants have similar
properties and can be extended to full invariants of virtual links.

\section*{Acknowledgement}
The author would like to thank Saziye Bayram for finding an error in a
calculation, and Dan Silver for informing him about his results.

\end{document}